\newcommand{\Z}{\mathbb{Z}}
\newcommand{\C}{\mathbb{C}}

\newcommand{\bea}{\begin{eqnarray}}
\newcommand{\eea}{\end{eqnarray}}
\newcommand{\beast}{\begin{eqnarray*}}
\newcommand{\eeast}{\end{eqnarray*}}

\newcommand{\E}{\mathbb{E}}

\documentclass[ECP]{ejpecp} 

\usepackage{xcolor}

\SHORTTITLE{Casimir Effect}
\TITLE{The Stochastic Casimir Effect} 

\AUTHORS{%
  Ruibo~Kou\footnote{Department of Mathematics, University of Warwick, UK.
    \EMAIL{Ruibo.Kou@warwick.ac.uk}}
  \and 
  Roger~Tribe\footnote{Department of Mathematics, University of Warwick, UK. \EMAIL{r.p.tribe@warwick.ac.uk}}
     \and 
  Oleg~Zaboronski\footnote{Department of Mathematics, University of Warwick, UK.
   \EMAIL{o.v.zaboronski@warwick.ac.uk} \orcid{0000-0001-7530-6643}}
    }

\KEYWORDS{Annihilating Brownian Motions, Casimir Effect} 

\AMSSUBJ{82C22, 81T55} 

\SUBMITTED{January 00, 0000} 
\ACCEPTED{January 00, 0000} 

\VOLUME{0}
\YEAR{2023}
\PAPERNUM{0}
\DOI{10.1214/YY-TN}

\ABSTRACT{We model the one-dimensional `classical' vacuum by a system of annihilating Brownian motions
on $\R$ with pairwise immigration. A pair of reflecting or absorbing walls placed in such a vacuum at separation
$L$ experiences an attractive force which decays exponentially with $L$. This phenomenon can be regarded as a purely 
classical Casimir
effect for a system of interacting Brownian motions. }

\newcommand{\R}{\mathbb{R}}

\begin{document}


\section{Introduction}\label{sec_intro}
The present note is motivated by a beautiful natural phenomenon arising due to
the complicated structure of the physical vacuum, the so called the Casimir effect predicted theoretically in \cite{Casimir predict} (see \cite{Casimir observe} for a review). 

Somewhat superficially,
the effect can be described as follows. According to the quantum field theory, the vacuum is teeming with virtual particles
popping in and out of existence. (Fundamentally, this is a consequence of the Heisenberg uncertainty principle between
energy $E$ and time $t$, $\Delta E\Delta t\geq \hbar/2$, which allows for existence of virtual particles which
violates energy conservation.) Imagine two walls placed in vacuum at separation $L$. The states of the virtual particles
are quantised in the region bounded by the walls and unconstrained outside of the walls. As a result, the density of the soup
of virtual particles in the bounded region is lower than outside and the walls experience an effective attraction which is
called the Casimir effect. Under the assumption that the virtual field is purely electromagnetic, it was shown in \cite{Casimir predict} that 
the strength of the attracting force decays as $L^{-4}$. The experimental verification of Casimir's theoretical prediction came half
a century later, see \cite{Casimir exper}.   

Our specific aim is to show that a Casimir-like effect occurs for a classical stochastic particle system. Our principal
example is the system of annihilating Brownian motions on $\R$ with
pairwise immigration, which can be regarded as the simplest `classical' model of the vacuum: here
particles are identical to anti-particles and it is assumed
that there are no carriers of interaction such as photons; the local creation of virtual particle-antiparticle pairs in the quantum case
is mimicked by a singular immigration mechanism - pairs of particles are immigrated to every point in space at an infinite rate to
prevent the immediate annihilation of all the pairs. The result is a soup of (identical) particles
which get created in pairs at a certain rate and annihilate in pairs on contact.  In between the collisions the particles
move independently of other particles and their motion is assumed to be Brownian.  
The walls are just the endpoints of an interval of length $L$ and both 
absorption and
reflection  will be considered as possible boundary behaviour. A steady state emerges where 
particle creation balances the particle annihilation and possible absorption.

We will calculate the expected force on the walls in the steady state for both the reflective and the absorbing boundary conditions  
and find that it is attractive and decays exponentially with $L$. 
The analysis is made possible due to the exact solvability of the interacting particle system at hand manifesting itself through
the existence of an infinite set of Markov dualities relating it to annihilating Brownian motions with finitely many particles. 
The rest of the paper is dedicated to the rigorous analysis of the model, but it is hard to resist mentioning that annihilating
Brownian motions have an experimental realisation as excitation in carbon nano-wires, see e.g. \cite{cnt} and
references therein. Might it be possible to
observe this `stochastic' Casimir effect in Nature? 

The paper is organised as follows. In section \ref{sec_mod} we describe the interacting particle system and
the boundary conditions associated with the walls, including the formulae
for the expected (`Casimir') force for both the reflective and the absorbing cases. The leading asymptotics  giving the dependence
of the expected Casimir force on the inter-wall separation are also given. The proofs can be found in section \ref{sec_prfs}.  
\section{The model and the main results.}\label{sec_mod}

We arrive at the simplest model of the vacuum - annihilating Brownian motions with pairwise immigration. 
More formally this can be constructed as a continuum limit of lattice models as follows.
Consider annihilating random walks on $\epsilon \Z$ with  immigration onto pairs of neighbouring sites. This is a continuous time Markov process 
$(\eta^{\epsilon}_t)_{t\geq 0}$ on $\Omega=\{0,1\}^{\epsilon \Z}$ fully characterised
by a deterministic initial condition and the following transition rates for the particles in pairs
of sites $\{x, x+\epsilon\}$: independently for each
$x\in \epsilon \Z$ we have transitions
\beast
(\eta (x), \eta(x+\epsilon))&\stackrel{1}{\rightarrow}& (0, \eta(x)\oplus \eta(x+\epsilon)) ~~\mbox{(Right jump and annihilation)},\\
(\eta (x), \eta(x+\epsilon))&\stackrel{1}{\rightarrow} &(\eta(x)\oplus \eta(x+\epsilon),0) ~~\mbox{(Left jump and annihilation)},\\
(\eta (x), \eta(x+\epsilon))&\stackrel{\beta/\epsilon^2}{\rightarrow}& (\eta(x)\oplus 1, \eta(x+\epsilon)\oplus 1) ~~\mbox{(Pairwise immigration)}.\\
\eeast
Here $\oplus$ denotes the addition modulo two, see e.g. \cite{Griffeath} for more details. Notice that the immigration rate approaches infinity in the limit 
$\epsilon \downarrow 0$, which insures existence of a non-trivial continuous limit. Define the following process
taking values in the space of counting measures on $\R$:
\beast
X^{\epsilon}_t(A)=\eta_{\epsilon^{-2}t}(A\cap \epsilon \Z), \quad A\in \frak{B}(\R),~t\geq 0,
\eeast
where $\eta(B):=\sum_{x\in B} \eta(x)$ for any $B\subset \epsilon \Z$. It was shown in \cite{ARWPI}
that the $\epsilon \rightarrow 0$ limit of any finite-dimensional distribution for the family of processes 
$(X^\epsilon_t)_{t\geq 0}$ exists. The resulting set of finite-dimensional distributions
is consistent and depends on $\beta$ in a non-trivial way.
We define annihilating Brownian motions with
pairwise immigration (ABMPI($\beta$)) as a process with this set of finite-dimensional dimensional distributions.
The present paper deals with the one-dimensional marginals only and
we do not discuss the regularity properties of the corresponding paths.
We also only consider the steady state of the process, so for definiteness we will only
consider ABMPI$(\beta)$ with {\it zero initial conditions}.

Consider ABMPI($\beta$) with either absorbing, or reflecting, walls placed at the points $0, L\in \R$. 
This means that particles approaching from either direction get absorbed (respectively reflected) at $0$ and $L$. As a result,
the particle subsystems contained in $(-\infty,0)$, $(0,L)$ and $(L, \infty)$ are independent for all $t>0$.
Such boundary behaviour can also
be defined as a limit of a particle system on $\epsilon \Z$ by adjusting the jumping rates the boundary points, see \cite{ARWPI} for details, where completely inhomogeneous systems are studied. 

Our aim is to build the simplest model of the Casimir force for this model. In the reflective case, it is natural
to assume that each reflected particle transfers a fixed amount of momentum to the wall. For the discrete model,
the number of collisions per unit of time is proportional to the particle number at the wall multiplied by the
rate of the jumps in the direction of the wall. The overall constant relating the particle number to the
force can be set to $1$ by a choice of measurement units. The above heuristics motivate the following definition, which depends on the function
$x\mapsto \rho(x)$ defined as the expected density of ABMPI($\beta$) particles in the steady state.
\begin{definition}\label{def of Casimir force R}(Casimir force at the reflective wall)
    For the system of particles that are reflected at $0$ and $L$, the expected Casimir force on the wall at $0$ is 
    \[F^R(L):=\rho(0-)-\rho(0+).\]
    \end{definition}
\noindent
By symmetry, the expected force exerted on the wall at $L$ is equal to $-F^R(L)$. Thus a positive value of 
$F^R(L)$ corresponds to attraction between the walls. 

Similarly, in the absorbing case each particle transfers momentum to the wall as it gets absorbed. Therefore
it is natural to model the force exerted on the wall by each of the particle subsystems by the 
flux of the absorbed particles into the wall. In turn the expected flux is proportional to the derivative of the density.
The situation in the continuum is complicated by the fact that this derivative is infinite at the wall.
This leads us to the following definition of the Casimir force in the absorbing case:
\begin{definition}(Casimir force at the absorbing wall)\label{defn42}
For the system of particles that are absorbed at $0$ and $L$, the expected Casimir force on the wall at $0$ is
    \[F^A(L):=\lim_{x\downarrow 0}\left(-\frac{d}{dx}\rho(-x)-\frac{d}{dx}\rho(x)\right).\]
\end{definition}

The goal of the present paper is to establish the Casimir effect for ABMPI$(\beta)$ by proving
that $F^R(L)\neq 0$, $F^A(L)\neq 0$.  Our model is so tractable that we can actually calculate 
$F^R(L)$ and $F^A(L)$ exactly, which allows also to describe the large $L$ behaviour.
Our main results can be stated as follows:
\begin{theorem}\label{Casimir force, reflected boundary}
    For ABMPI($\beta)$, where particles are reflected at $0$ and $L$, the expected Casimir force 
    is given by
    \begin{equation}
        F^{R}(L)=\sqrt{\frac{\beta}{2}}\frac{1-e^{-L\sqrt{2\beta}}}{\sinh L\sqrt{2\beta}} \geq 0.
    \end{equation}
    Moreover, ${F^R(L)}\approx {e^{-L\sqrt{2\beta}}}$ as $L\to \infty.$
\end{theorem}
Here we used the notation "$\approx$" to describe the weak asymptotic relation between  functions, that is, for two given functions $f(x)$ and $g(x)$, we say $f(x)\approx g(x)$ if and only if
\[\lim_{x\to\infty}\frac{f(x)}{g(x)}=C>0.\]

\begin{theorem}\label{Casimir force, absorbed boundary}
    For ABMPI($\beta$), where particles are absorbed at $0$ and $L$, the expected Casimir force is given by
       \begin{equation}\label{finaly express for F_L^A}
       \begin{split}
                   F^{A}(L)=\frac{2\beta}{\pi}\int_{0}^\infty\frac{e^{-4\beta L^2y}}{y}\left(1-\theta^2\left(0,\frac{i}{\pi y}\right)\right)dy+2\beta\int_{0}^\infty e^{-4\beta L^2y}\theta^2\left(\frac{1}{2},i\pi y\right)dy,
       \end{split}
    \end{equation}
    where $\theta:\C\times \mathbb{H}\rightarrow \C$ is the theta-function.
    Moreover, $F^A(L)\approx\frac{1}{\sqrt{L}}e^{-L\sqrt{8\beta}},$ as $L\to\infty$.
\end{theorem}
Notice that in each of the cases considered the expected Casimir force appears to be attractive, see Figs. 1,2 for an illustration
of the $\beta=1$ case. It would be interesting to study the full law of the Casimir force as a process in time. 
Such an investigation should be facilitated by the fact that ABMPI($\beta$) is an extended Pfaffian point process,
see \cite{eppp}, but at present we restrict ourselves to the study of the first moment of the force at $t=\infty$ only.

The rest of the paper is dedicated to the proofs of the stated theorems.
\begin{figure}[htbp]
\centering 
\begin{minipage}[t]{0.48\textwidth}
    \centering
        \includegraphics[width=6cm]{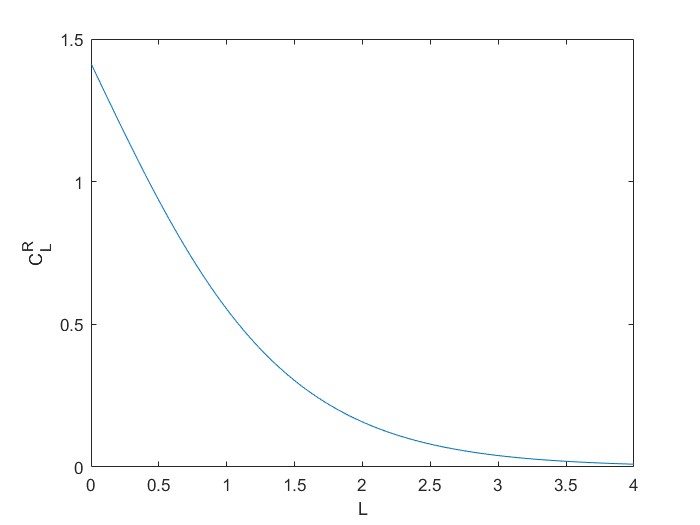}
        \caption{Casimir force, reflected case}
\end{minipage}
\begin{minipage}[t]{0.48\textwidth}
    \centering
        \includegraphics[width=6cm]{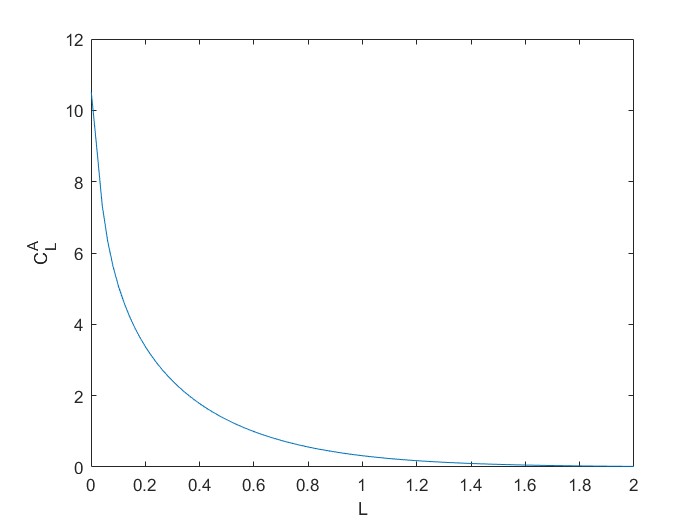}
        \caption{Casimir force, absorbed case}
\end{minipage}

\label{fi:myfigure}
\end{figure}

\section{Proofs}\label{sec_prfs}
\subsection{Theorem \ref{Casimir force, reflected boundary}}
 Let $V_t(x):=\E\left[(-1)^{X_{t}(0,x)}\right]$ be the expected parity of particles in the interval $(0,x)$
 at time $t>0$. The parity of the number of particles in a Borel subset of $\R$ is a Markov duality function
 which can be used to prove that the one-dimensional law of ABMPI$(\beta)$ 
 is a Pfaffian point process for any deterministic initial condition, including the empty initial condition
 at hand, see \cite{benam}, \cite{tz1}, \cite{ARWPI} for details. The knowledge of the expected
 parity in an interval allows one to calculate the expected density:
 \bea\label{vtorho}
 \rho_t(0+)=-\frac{1}{2}\partial_xV_t(0+).
 \eea
 In particular, $\rho(0+)=-\frac{1}{2}\partial_xV(0+)$, where $V(x)=\lim_{t\rightarrow \infty} V_t(x)$
 is the expected parity of the number of particles in $(0,x)$ in the steady state.
 
 The derivation of (\ref{vtorho}) requires some care as the function under the sign of expectation
 is not differentiable, see \cite{tz1} for details. It was shown in \cite{ARWPI}, by passing to the limit in lattice equations, that the function $(x,t)\mapsto V_t(x)$
 satisfies the following equation
 \bea \label{pde10}
 \begin{cases}
 \left(\partial_t- \partial_x^2+2\beta \right)V_t(x)=0,& 0<x<L,t>0,\label{veqn}\\
 V_t(0+)=1,& t>0,\label{vlbc}\\
 V_0(x)=1, &0<x<L,\label{vbc}
 \end{cases}
 \eea
which expresses the fact that $(-1)^{X_t(0,x)}$ is a Markov duality function between ABMPI($\beta$) and the system of
two annihilating particles started at $(0,t)$ and $(x,t)$ and run backwards in time.
 We are not going to repeat the derivation of (\ref{veqn}) here, just notice that the $\beta=0$ case in the continuum has been
 treated in \cite{tz1}. The appearance of the $\beta$-term is easy to see before taking the continuous limit
 by noticing that only the immigration of a pair straddling the boundary can change the parity of the number of particles
 in the interval. The boundary condition (\ref{vlbc}) follows from the fact that $X_t(0,0)=0$ a.s. and the dominated
 convergence theorem. Finally, the initial condition (\ref{vbc}) follows from the fact that we start our process from the
 empty state. 
 
 To finish formulating the boundary value problem characterising $V$ completely one needs a boundary
 condition at $x=L$. This can be done using the following conservation law: 
the evolution of ABWPI($\beta$) with reflective boundary conditons at $0,L$ preserves the parity of $X_t(0,L)$ as particles in the interval
$(0,L)$ can only appear and disappear in pairs. Therefore, $(-1)^{X_t (0,L)}=(-1)^{X_0(0,L)}=1$ a.s. 
The p.d.e. (\ref{pde10}) has a limit as $t \to \infty$, which can be used to establish that the process converges to a steady state.
Using this observation together with (\ref{veqn}), (\ref{vlbc})
we conclude that the expected parity in the steady state $V(x)$ is uniquely determined by 
\bea \label{fineqrefl}
\begin{cases}
 \left(\partial_x^2-2\beta \right)V(x)=0,& ~0<x<L,t>0,\\
 V(L-)=V(0+)=1,& ~t>0.
\end{cases}
\eea 
The solution is 
\beast
V(x)=\frac{\sinh \sqrt{2\beta}x+\sinh\sqrt{2\beta}(L-x)}{\sinh\sqrt{2\beta}L}, 0<x<1,
\eeast
substituting which into (\ref{vtorho}) one finds
\bea\label{rhoreflL}
\rho(0+)=\sqrt{\frac{\beta}{2}}\frac{\cosh \sqrt{2\beta}L-1}{\sinh \sqrt{2\beta L}}.
\eea
To find the density $\rho(0-)$ just outside of the wall one can either solve (\ref{fineqrefl}) on $(0,\infty)$ under the condition
that the solution is bounded at infinity, or simpler still, take the $L\rightarrow \infty$ limit in (\ref{rhoreflL}). Using either method, one finds
\bea\label{rhoreflinf}
\rho(0-)=\sqrt{\frac{\beta}{2}}.
\eea
Using the definition (\ref{def of Casimir force R}) one obtains the final answer for the expected Casimir force
\beast
F^R(L)=\sqrt{\frac{\beta}{2}}\frac{1-\exp\left(- \sqrt{2\beta}L\right)}{\sinh \sqrt{2\beta L}}.
\eeast
Clearly, $F^R(L)\approx e^{-\sqrt{2\beta}L}$ for $L\rightarrow \infty$.
Theorem \ref{Casimir force, reflected boundary} is proved.
\subsection{Theorem \ref{Casimir force, absorbed boundary}}
The calculation of the Casimir force for the absorbing case is more complicated due to the absence of parity
conservation and the need to know $\rho(x)$ for $x\neq 0$, see the definition \ref{defn42}.

To generalise the discussion of the previous section, let
$V_t(x,y):=\E\left[(-1)^{X_{t}(x,y)}\right]$ and $V(x,y):=\lim_{t\rightarrow \infty}V_t(x,y)$
be the parity of the number of particles in $(x,y)\in \R$ at time $t>0$ and in the steady state correspondingly. 
Then
\bea\label{vtorhox}
\rho_t(x)=-\frac{1}{2}\partial_y V_t(x,y)\mid_{y=x},
\eea
which is a counterpart of (\ref{vtorho}).
For $0\leq x\leq y\leq L$ and $t=\infty$, $V(x,y)$
is characterised as the unique solution to the following boundary value problem:
\begin{equation}\label{bvpabs}
\begin{cases}
\left(\partial_x^2+\partial_y^2 \right) V(x,y)-4\beta V(x,y)=0,\quad &\text{for}\quad 0<x<y<L,\\
V(x,x)=1\quad &\text{for}\quad 0<x< L,\\
\partial_x V(x,y)\big|_{x=0}=0\quad &\text{for}\quad 0< y< L,\\
\partial_y V(x,y)\big|_{y=L}=0\quad &\text{for}\quad 0< x< L.
\end{cases}
\end{equation}
The equation for $V$ in the bulk along with the boundary condition 
on the diagonal $x=y$
have been derived in \cite{ARWPI}, see also a  brief discussion
in the previous section.  
The boundary conditions at $x=0$ and $y=L$ require a separate discussion. By symmetry, it is 
sufficient to consider the case $x=0$. Consider annihilating random walks with pairwise immigration
on $\epsilon \Z_{+}$ as discussed in the Introduction. An absorbing wall at $0$ corresponds to the following 
hopping rates: a particle at zero can jump to $\epsilon$ at rate $1$ and killed at rate $1$. Killing can be represented
as a hop to a virtual site $-\epsilon$, the hopping rates out of which are equal to zero. A calculation using the corresponding
infinitesimal generator gives
\beast
\partial_t V_t(0,y)=V_t(\epsilon,y)-V_t(0,y)+\ldots,
\eeast
where $\ldots$ denotes the contribution of jumps and immigration at $y$. Using the virtual site $-\epsilon$ this
be rewritten as $\partial_t V_t(0,y)=(V_t(-\epsilon,y)+V_t(\epsilon,y)-2V_t(0,y))+\ldots$ provided
\beast
V(0,y)-V(-\epsilon,y)=0.
\eeast
In the continuous limit $\epsilon \rightarrow 0$ this gives the  Neumann boundary condition at $0$,
whereas the remaining operator acting at $0+$ becomes the Laplacian, exactly as in the bulk. 
This consideration is an example of the argument replacing the inhomogeneity in the hopping rates by an effective boundary
condition which is the first step in constructing the coordinate Bethe Ansatz, see e.g. \cite{borodin} for previous applications to
interacting particle systems. We remark that the boundary condition at hand also makes sense from the point of duality: the system of annihilating
Brownian motions in $[0,L]$ \emph{absorbed} at the boundary is dual to the system of annihilating
Brownian motions in reverse time \emph{reflected} at the boundary and the Neumann boundary condition
corresponds to reflection. 

Similarly to the reflective case, the expected parity of the number of particles 
in the subinterval $(x,y)$ of $(-\infty,0)$ in the steady state can be characterised
as the bounded solution of 
 \begin{equation}\label{bvpabsinf}
\begin{cases}
\left(\partial_x^2+\partial_y^2\right) V(x,y)-4\beta V(x,y)=0,\quad &\text{for}\quad -\infty<x<y<0,\\
V(x,x)=1\quad &\text{for}\quad x< 0,\\
\partial_y V(x,y)\big|_{y=0}=0\quad &\text{for}\quad x<0.
\end{cases}
\end{equation}

The solution to (\ref{bvpabsinf}) is

\begin{equation}\label{vout}
V(x,y)=1+\int_{0}^\infty4\beta e^{-4\beta u}\text{erf}\left(\frac{x+y}{2\sqrt{2u}}\right)
\text{erf}\left(\frac{y-x}{2\sqrt{2u}}\right)du,
\end{equation}
where $x\mapsto \text{erf}(x)=\frac{2}{\sqrt{\pi}}\int_0^x e^{-u^2}du$ is the error function. 
Checking that (\ref{vout}) satisfies the boundary conditions in (\ref{bvpabsinf}) is elementary, 
the check of the equation is also straightforward using the relation 
$$(\partial_u-2\partial_x^2)\text{erf}\left(\frac{x}{2\sqrt{2u}} \right)=0.$$
Substituting (\ref{vout}) into (\ref{vtorhox}) one gets the steady-state particle density in half-space:
\bea\label{ehossa}
\rho(x)=-\int_{0}^\infty  \frac{2\beta e^{-4\beta u}}{\sqrt{2\pi u}}\text{erf}\left(\frac{x}{\sqrt{2u}}\right)du,~x<0.
\eea

Next, let us find the expected particle density in $[0,L]$ by solving (\ref{bvpabs}). 
Notice that the function $U:=V-1$ satisfies the same boundary conditions 
as $V$ at $x=0$ and $y=L$, and vanishes on the diagonal $x=y$. It satisfies the equation
$\left(\partial_x^2+\partial_y^2 -4\beta\right) U(x,y)=4\beta $ in the wedge
$0<x<y<L$, which can be solved by a skew-symmetric extension to the domain
$[0,L]^2$. Namely, $U$ can be characterised as the restriction to the wedge of the
solution to the following problem: 
\begin{equation}\label{bvpabsext}
\begin{cases}
\left(\partial_x^2+\partial_y^2-4\beta \right) V(x,y)=4\beta \text{sign}(y-x),\quad &\text{for}\quad 0<x,y<L,\\
\partial_x V(x,y)\big|_{x=0,L}=0\quad &\text{for}\quad 0< y< L,\\
\partial_y V(x,y)\big|_{y=0,L}=0\quad &\text{for}\quad 0< x< L.
\end{cases}
\end{equation}
The equation (\ref{bvpabsext}) can be solved by expanding the solution w.r.t. the complete
system $\left(\cos \left(\frac{\pi mx }{L}\right)\cos\left(\frac{\pi ny}{L}\right)\right)_{m,n\geq 0}$
of the eigenfunction of the Laplacian on the $[0,L]^2$ with Neumann boundary conditions. 
The answer is
\begin{equation*}
        V(x,y)=1+\sum_{m,n\geq 0} a_{m,n}\cos{\frac{m\pi x}{L}}\cos{\frac{n\pi y}{L}},~0<x<y<L
        \end{equation*}
where $a_{m,n}=-a_{n,m}$ are the Fourier coefficients given by
\begin{equation*}
    a_{m,n}=
    \begin{cases}
    \frac{-16\beta L^2(1-\cos{(m\pi)})}{m^2\pi^2(m^2\pi^2+4\beta L^2)},\quad&\text{for}\quad m>0,n=0,\\
    \frac{-16\beta L^2(\cos{(n\pi)}-1)}{n^2\pi^2(n^2\pi^2+4\beta L^2)},\quad&\text{for}\quad m=0,n>0,\\
    \frac{-32\beta L^2(\cos{((m+n)\pi)}-1)}{(n^2-m^2)\pi^2((m^2+n^2)\pi^2+4\beta L^2)},\quad&\text{for}\quad m,n>0,m\neq n,\\
    0\quad&\text{for}\quad m=n.
    \end{cases}
\end{equation*}
Using (\ref{vtorhox}) one can conclude that the expected particle density  is
\begin{equation*}
\begin{split}
        \rho(x)&=\sum_{n=1}^\infty\frac{8\beta L(1-\cos{(n\pi)}}{n\pi(n^2\pi^2+4\beta L^2)}\sin{\frac{n \pi x}{L}}\\
        &+\sum_{m,n> 0}\frac{16\beta L n(1-\cos{((m+n)\pi)})}{\pi(n^2-m^2)((m^2+n^2)\pi^2+4\beta L^2)}\cos{\frac{m\pi x}{L}}\sin{\frac{n\pi x}{L}},0<x<L.
\end{split}
\end{equation*}
The term-by-term differentiation of the Fourier series is justified as $V$ is a $C^1$ function. 

Using  the definition \ref{defn42} of the Casimir force applied to the wall at zero one finds
that 
\bea\label{finst1}
F^A(L)=\lim_{x\downarrow 0} J(-x)-J(x),
\eea
where the flux $J$ is
\bea\label{finstflux}
J(x)=
\begin{cases}
\int_0^\infty\frac{2\beta }{\pi u}e^{-\frac{x^2}{2u}}e^{-4\beta u}du,&\!\!\!\!\text{for}~ x<0,\\
\sum\limits_{n>0}\frac{8\beta(1-\cos n\pi)}{n^2\pi^2+4\beta L^2}\cos \frac{n\pi x}{L}+\!\!\!\!\sum\limits_{m,n>0}
\!\!\!\!\frac{8\beta(1-\cos(m+n)\pi)}{(m^2+n^2)\pi^2
+4\beta L^2}\cos\frac{m\pi x}{L}\cos\frac{n\pi x}{L},&\!\!\!\!\text{for}~ x\in(0,L).
\end{cases}
\eea
The expression for the flux for $x\in (0,L)$ can be written in a more compact form by firstly observing
that 
\[
\frac{1}{(m^2+n^2)\pi^2+4\beta L^2}=\int_0^\infty e^{-((m^2+n^2)\pi^2+4\beta L^2)u}du, ~L>0,
\] 
and secondly using that
\[
\sum_{n=1}^\infty e^{-n^2\pi^2 \frac{u}{L^2}}\cos \frac{n\pi x}{L}=\frac{1}{2}
\left(\theta\left(\frac{x}{2L},\frac{i\pi u}{L^2} \right)-1 \right),
\]
where $(z,\tau)\mapsto \theta (z,\tau)$ is the Jacobi theta function defined on $\C\times \mathbb{H}$, see \cite{theta},
\cite{E1 function} for review. The result is 
\bea\label{finst2}
J(x)=\frac{2\beta}{L^2}\int_0^\infty e^{-4\beta u} \left(\theta^2
\left(\frac{x}{2L},\frac{i\pi u}{L^2}\right)-\theta^2
\left(\frac{1}{2}, \frac{i\pi u}{L^2} \right) \right)+o(x^0),~0<x<L.
\eea
In the analysis presented below we will use the following representation
for the theta function:
\bea\label{finst5}
\theta\left(0,\frac{i}{\pi u} \right)=1+2e^{-\frac{1}{u}}+2e^{-\frac{2}{u}}R(u),
\eea
where $R(u)=\sum_{n\geq 2}e^{-\frac{n^2}{u}}$. As it is easy to check,
\beast
0\leq R(u)\leq \sum_{n\geq 0}e^{-\frac{4n}{u}}=\frac{1}{1-\exp\left(-\frac{4}{u}\right)}, ~u>0.
\eeast
Moreover, 
\beast
R(u)=\left(1+\frac{u}{4}\right)B(u),
\eeast
where the function $B$ is bounded on $(0,\infty)$ by some constant $M>0$.

To facilitate the calculation of the limit in (\ref{finst1}), one can rewrite $\theta^2
\left(\frac{x}{2L},\frac{i\pi u}{L^2}\right)$ using the Jacobi identity:
\bea\label{finst3}
\theta^2
\left(\frac{x}{2L},\frac{i\pi u}{L^2}\right)=\frac{L^2}{\pi u} e^{-\frac{x^2}{2u}} \theta^2
\left(\frac{xL}{2\pi i u}, \frac{iL^2}{\pi u}\right).
\eea
Substituting (\ref{finst3}) into (\ref{finst2}), and then substituting the result along with the first line
 of (\ref{finstflux}) into (\ref{finst1}), one finds that
 \beast
 F^A(L)=\lim_{x\downarrow 0}\frac{2\beta}{\pi}\int_0^\infty \frac{e^{-4\beta u}}{u}e^{-\frac{x^2}{2u}}
 \left[1-\theta^2\left(\frac{xL}{2\pi i u}, \frac{iL^2}{\pi u}\right)\right]du+\frac{2\beta}{L^2}\int_0^\infty e^{-4\beta u}\theta^2(\frac{1}{2},\frac{i\pi u}{L^2})du,
 \eeast
 Noticing with the help of (\ref{finst5}) that 
 $$
 e^{-\frac{x^2}{2u}}\left| 1-\theta^2\left(\frac{xL}{2\pi i u}, \frac{iL^2}{\pi u}\right)\right|
 \leq g(u):=4\frac{e^{-\frac{L^2}{2u}}}{\left(1-e^{-\frac{L^2}{2u}}\right)^2},~x\in (0,L/2),u>0,
$$
and that $\int_0^\infty \frac{e^{-4\beta u}}{u}g(u)du<\infty$, one can use the dominated convergence
theorem to justify the interchange of the $\lim_{x\downarrow 0}$ and the integration
over $u$
in the last expression for $F^A(L)$.
The
result is
\bea\label{finst4}
F^A(L)=\frac{2\beta}{\pi}\int_0^\infty \frac{e^{-4L^2\beta u}}{u}
 \left[1-\theta^2\left(0, \frac{i}{\pi u}\right)\right]du+2\beta\int_0^\infty e^{-4L^2\beta u}\theta^2\left(\frac{1}{2},i\pi u\right)du,
\eea
which proves the statement 
(\ref{finaly express for F_L^A}) of the Theorem \ref{Casimir force, absorbed boundary}.

It remains to establish the large-$L$ asymptotic of (\ref{finst4}). 
Fix any $\delta>0$. The representation (\ref{finst5})
leads to the estimate
\bea\label{finst6}
\left| \int_\delta^\infty \frac{e^{-4L^2\beta u}}{u}
 \left[1-\theta^2\left(0, \frac{i}{\pi u}\right)\right]du\right|\leq \frac{C_1(M,\delta)}{L^2}e^{-4L^2\beta \delta},
\eea
for some positive $L$-independent constant $C_1(M,\delta)$. Similarly, there is a positive constant $C_2(\delta)$
such that
\bea
\left| \int_0^\delta \frac{e^{-4L^2\beta u}}{u}
 \left[1-\theta^2\left(0, \frac{i}{\pi u}\right)\right]du\right|\leq C_2(\delta)\int_0^\delta \frac{e^{-4L^2\beta u-\frac{1}{u}}}{u}du\nonumber\\
\leq   C_2(\delta)\int_0^\infty \frac{e^{-L\left(4\beta u+\frac{1}{u}\right)}}{u}du=\frac{C_2(\delta)}{\sqrt{8\pi \beta^{1/2}L}}
e^{-4\sqrt{\beta}L}\left(1+O(L^{-1/2}) \right).\label{finst7}
\eea
The last integral in the above computation is evaluated using the Laplace method (see \cite{Laplace} for a review) by noticing that the exponent
$F(u):=4\beta u+1/u$  achieves a global non-degenerate minimum on $(0,\infty)$ at $u_c=1/\sqrt{4\beta}$.  
The inequalities (\ref{finst6}, \ref{finst7}) are sufficient for the estimation of the first integral in the r.h.s. of (\ref{finst4}).
The estimate of the second integral is done in a similar manner provided that the Jacobi identity is used to rewrite
the theta function in the form suitable for the application of the Laplace method:
\beast
\theta^2\left(\frac{1}{2},i\pi u \right)=\frac{e^{-\frac{1}{2u}}}{\pi u} \theta^2\left(\frac{i}{2\pi y},\frac{i}{\pi u} \right).
\eeast 
Observing that $\lim_{u\downarrow 0} \theta\left(\frac{i}{2\pi y},\frac{i}{\pi u} \right)=2$, one finds eventually that
\bea\label{finst8}
\int_0^\infty e^{-4L^2\beta u}\theta^2\left(\frac{1}{2},i\pi u\right)du\approx\!\! 
\int_0^\infty \frac{e^{-4L^2\beta u-\frac{1}{2u}}}{u}du
=\!\!\int_0^\infty \frac{e^{-L\left(4\beta u+\frac{1}{2u}\right)}}{u}du
\approx \!\!\frac{e^{-L\sqrt{8\beta}}}{\sqrt{L}},
\eea
where the Laplace method is used in the last step. Comparing (\ref{finst6}), (\ref{finst7}) and (\ref{finst8})
one finds that
\beast
F^{A}(L)\approx \frac{e^{-\sqrt{8\beta}L}}{\sqrt{L}}.
\eeast 
Theorem \ref{Casimir force, absorbed boundary} is proved. 

\begin{acks}
We are grateful to Alexandre Povolotsky for many illuminating discussions.
\end{acks}


\end{document}